%% file: main.tex
\documentclass{article}

\usepackage[nonatbib, preprint]{neurips_2022}

\usepackage[utf8]{inputenc} 
\usepackage[T1]{fontenc}    
\usepackage{hyperref}       
\usepackage{url}            
\usepackage{booktabs}       
\usepackage{amsfonts}       
\usepackage{nicefrac}       
\usepackage{microtype}      
\usepackage{xcolor}         
\usepackage{amsmath}
\usepackage{tikz}
\usetikzlibrary{positioning}
\usepackage{amssymb}
\usepackage[noabbrev, capitalise]{cleveref}
\usepackage{graphicx}
\usepackage{float}
\usepackage{subfig}
\usepackage{amsthm,thmtools}
\usepackage{wrapfig}

\usepackage[
backend=biber,
style=alphabetic,
sorting=ynt
]{biblatex}

\newtheorem{theorem}{Theorem}
\input{Tikz_pictures/tikz}

\title{Spectral Neural Operators}

\author{%
    Vladimir Fanaskov\\
    Center for Artificial Intelligence Technology\\
    Skolkovo Institute of Science and Technology\\
    Bolshoy Boulevard 30, bld. 1 Moscow, Russia, 121205\\
    \texttt{v.fanaskov@skoltech.ru}\\
    \And
    Ivan Oseledets\\
    Artificial Intelligence Research Institute (AIRI)\\
    \\
    Center for Artificial Intelligence Technology\\
    Skolkovo Institute of Science and Technology\\
    Bolshoy Boulevard 30, bld. 1 Moscow, Russia, 121205\\
    \texttt{i.oseledets@skoltech.ru}
}

\begin{document}

\maketitle

\begin{abstract}
    A plentitude of applications in scientific computing requires the approximation of mappings between Banach spaces. Recently introduced Fourier Neural Operator (FNO) and Deep Operator Network (DeepONet) can provide this functionality. For both of these neural operators, the input function is sampled on a given grid (uniform for FNO), and the output function is parametrized by a neural network. We argue that this parametrization leads to 1) opaque output that is hard to analyze and 2) systematic bias caused by aliasing errors in the case of FNO. The alternative, advocated in this article, is to use Chebyshev and Fourier series for both domain and codomain. The resulting Spectral Neural Operator (SNO) has transparent output, never suffers from aliasing, and may include many exact (lossless) operations on functions. The functionality is based on well-developed fast, and stable algorithms from spectral methods. The implementation requires only standard numerical linear algebra. Our benchmarks show that for many operators, SNO is superior to FNO and DeepONet.
\end{abstract}

\section{Introduction}
\label{section:Introduction}
Fourier Neural Operator \cite{li2020fourier} and Deep Operator Network \cite{lu2019deeponet} are two of the most well-developed neural operators to date. Benchmarks indicate that PDE solvers based on them can be superior to classical techniques at least in the regimes of low accuracy \cite{wang2021learning}, \cite{kovachki2021neural}, \cite{grady2022towards}, \cite{pathak2022fourcastnet}, \cite{yang2022generic}. Architectures of FNO and DeepONet are quite different, but they both use sampling to approximate input function and neural networks to approximate the output. We have reasons to believe that these choices are unfortunate when PDEs are considered.

\begin{figure*}
    \centering
    \subfloat[]{
        \label{fig:bad_network}
        \includegraphics[width=0.45\textwidth]{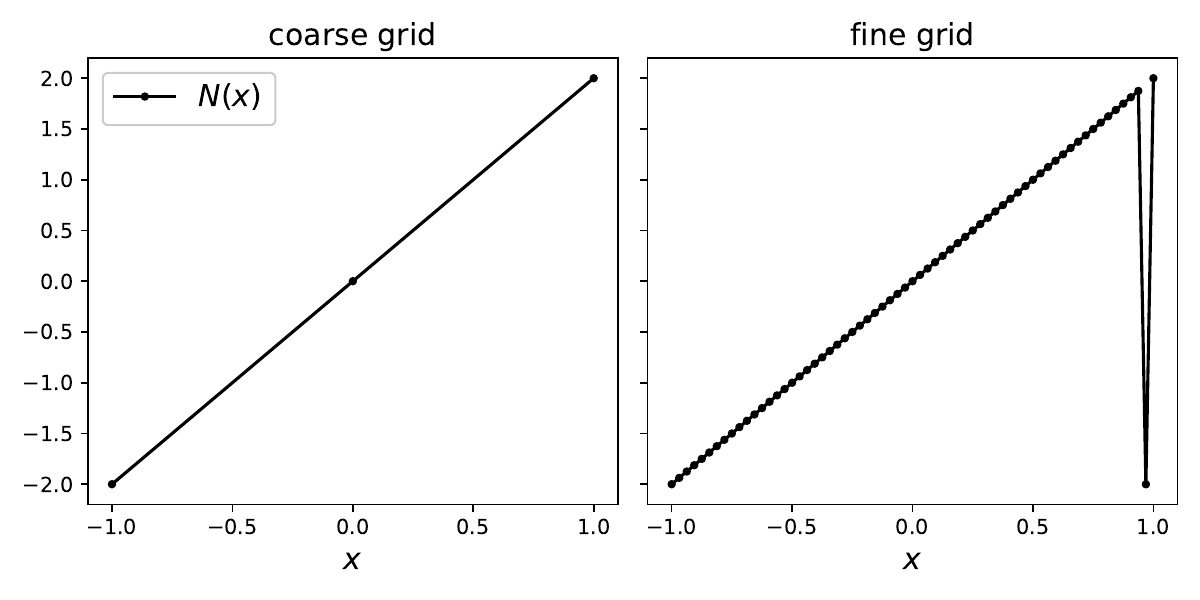}
    }
    \subfloat[]{
        \label{fig:inconsistent_FNO}
        \includegraphics[width=0.51\textwidth]{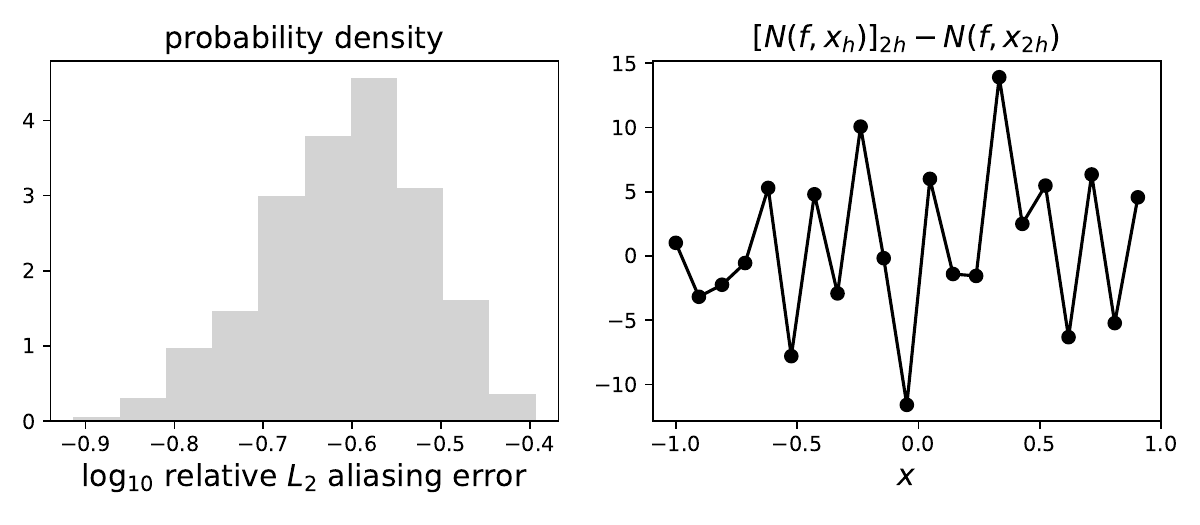}
    }
    \caption{(a) the output of neural network $N(x)$ computed on coarse and fine grids. On each subgrid, loss and gradients are zero, so the network provides the best (alas, pathological) approximation to $f(x) = 2x$ on the interval $[-1, 1]$ for all but the finest grid. (b) FNO trained to take the first derivative of random trigonometric polynomials on grid $2h$ is evaluated for the same functions on grid $h$. The difference between outputs for a particular function is on the right; statistics for $750$ functions is on the left. The average discrepancy is about $25\%$ of $L_2$ norm. Inconsistency of FNO results from aliasing of high harmonics produced by activation functions.}
    \label{fig:sno_motivation}
\end{figure*}

First, the output of the neural operator is a neural network, i.e., essentially a black-box function. In classical numerical methods such as FEM \cite{ciarlet2002finite}, spectral methods \cite{boyd2001chebyshev} and others, the parametrization allows for extracting bounds on function,  its derivatives, and any other local or global information in a controlled transparent manner. For neural networks, this is not the case. Consider neural networks with $\text{ReLU}$ activation functions in $1D$. The network is a piecewise linear function, so on a given interval it is uniquely fixed by partition and the value of derivative on each subinterval \cite{devore2021neural}. The number of subintervals is exponential in the number of layers, so the partition is not easy to find. Clearly, simple strategies such as sampling are not going to work. The example of this situation is given on \cref{fig:bad_network}. Here on grids with $h$, $h/2$, $h/4$, $\dots$, neural network approximates $2x$, and on the grid with $h/2^{n}$, it significantly diverges from $2x$. Such behavior is clearly undesirable from the perspective of scientific computing.\footnote{Perhaps the only argument for neural networks is to use them for high-dimensional problems \cite{yu2018deep}, \cite{grohs2018proof}, \cite{sirignano2018dgm}. However, in this context, other well-known parametrizations are available \cite{bungartz2004sparse}. Besides that, current neural operators are inapplicable to high-dimensional problems.}

Second, the representation of the input function by sampling implicitly restrict the class of functions and possible operations with them. Indeed, it is known that for periodic functions sampling on grid with $2N+1$ points allows to resolve frequencies up to $N$ (Nyquist frequency).\footnote{For non-periodic functions the same argument can be made using Chebyshev polynomials.} All higher frequencies with $k>N$  are indistinguishable from lower frequencies with $k\leq N$ (this is known as aliasing \cite[Chapter 11, Definition 23]{boyd2001chebyshev}). For DeepONet this implies it is defined on trigonometric polynomials with degree $\leq N$. For FNO situation is more difficult, because even when the input function is fully resolved, nonlinear activations are going to introduce higher frequencies, so the lower frequencies will be distorted. As a result FNO $N(u(x), x)$ computed on grid $x_{2h}$ produces different result when computed on grid $x_{2h}$. This is illustrated in \cref{fig:inconsistent_FNO}, where the aliasing error $\left\|N(u(x_{2h}), x_{2h}) - \left[N(u(x_{h}), x_{h})\right]_{2h}\right\|$ ($\left[\cdot\right]_{2h}$ is a projection on grid $2h$) is computed.

Our solution to both of outlined problems is to detach the process of function approximation (interpolation) from the mapping between functional spaces, and explicitly fix the highest possible resolution. To do that, for both domain and codomain of neural operator we consider functions represented by finite series of the form $\sum c_n f_n(x)$ where $f_n$ are either Chebyshev or trigonometric polynomials. Now, all functions are equivalent to vectors of finite length, and we can realize mapping between functional spaces with (almost) vanilla feedforward neural network. Benefits of this strategy include: 1) reach set of operations (integration, differentiation, shift) available for Chebyshev and trigonometric polynomials; 2) lossless operations on functions; 3) interpolation on the fine mesh can be performed with arbitrary accuracy; 4) availability of global information about obtained functions (f.e. bounds on functions and their derivatives); 5) well-studied approximation properties; 6) good compression for smooth functions; 7) overall similarity with spectral and pseudospectral methods.

Main contributions are:
\begin{itemize}
    \item Analysis of aliasing error in FNO (\cref{subsection:Activation_functions_and_aliasing}, \cref{appendix:aliasing_and_relu}).
    \item Definition of super-resolution and generalization tests (\cref{subsection:Superresolution}, \cref{appendix:Additional_comments_on_super_resolution}).
    \item New type of Spectral Neural Operators (SNO), based on Chebyshev and Fourier series (\cref{section:Spectral_neural_operator}, \cref{appendix:details_on_training_and_architectures}).
    \item Comparison of SNO, DeepONet and FNO on $16$ problems (\cref{section:Numerical_examples}, \cref{appendix:benchmarks_for_additional_architectures}).
\end{itemize}

\section{Chebyshev and trigonometric polynomials}
\label{section:Chebyshev_and_trigonometric_polynomials}
As we stated, we are going to consider neural networks operating on functions of the form $f(x) = \sum c_n f_n(x)$, $x\in[-1, 1]$\footnote{As usual, basis functions  for higher dimensions are constructed as direct product of $D=1$ basis functions. From now on we consider $D=1$ having in mind that generalization to higher dimensions is straightforward.}, where $f_n(x)$ are either Chebyshev polynomials defined by reccurence \cite[Section 1.2.1]{mason2002chebyshev} $T_{n}(x) = 2x T_{n-1}(x) - T_{n-2}(x)$, $T_{0}(x)=1$, $T_{1}(x)=x$ or trigonometric polynomials $e^{i\pi kx}$, $k\in \mathbb{Z}$. This particular choice of basis warrants some explanation, which we provide in the present section.

First, coefficients $c_n$ can be efficiently approximated by sampling. The standard way to obtain coefficients is to use orthogonality of $f_n$ with respect to $L_2$ norm and compute the integral $\int f_{n}(x)f(x)dx$. For Chebyshev polynomials we can sidestep the integration using samples on Chebyshev grid $x_{k} = \cos(k\pi \big/ n),~k=0,\dots,n$, trigonometric identity $T_{n}(\cos(\theta)) = \cos(n\theta)$ and Discrete Cosine Transform  (DCT-I) \cite[Section 4.7]{mason2002chebyshev}. For trigonometric polynomials, coefficients can be obtained in a similar way but with sampling on the uniform grid and Fast Fourier Transform (FFT) \cite[Chapter 3]{trefethen2000spectral}. So the cost to recover $N$ coefficients is $O(N\log(N))$ plus the cost of $N$ evaluations of $f(x)$. The procedure also come with reliable adaptive algorithms which can chop the series based on accuracy requirements \cite{aurentz2017chopping}, \cite[Chapter 3]{boyd2014solving}, \cite{trefethen2007computing}. Moreover, the very same DCT-I and FFT can be used to perform interpolation both on Chebyshev and the uniform grid with the same complexity $O(N\log(N))$ \cite{boyd1992fast}. Altogether, we have a duality, between values on the grid and values of coefficients. This duality is significantly exploited in pseudospectral methods, when one needs to apply nonlinear transformation to given Chebyshev or trigonometric series.

The second reason is efficient compression of smooth functions. It is known that for functions with $m$ continuous derivatives $n$-th coefficient is $O(n^{-m})$ for both Chebyshev \cite[Theorem 5.14]{mason2002chebyshev} and Fourier series \cite[Chapter 2, Miscellanious theorems and examples]{zygmund2002trigonometric}.

Besides that, many operations such as integration, differentiation, shift (for Fourier series), multiplication of two series can be performed with machine accuracy \cite[Equation (2.49)]{mason2002chebyshev}, \cite[Equation (2.43)]{mason2002chebyshev}, \cite[Chapter 3]{trefethen2000spectral}, \cite{giorgi2011polynomial}.

All these properties combined make trigonometric and Chebyshev polynomials especially suitable for PDE solutions by spectral methods \cite{boyd2001chebyshev} and adaptive lossless computations with functions \cite{trefethen2007computing}. The later goal is fully realized in the Chebfun software.\footnote{\href{https://www.chebfun.org}{https://www.chebfun.org}} Chebfun demonstrates that computations with continuous functions can be performed with essentially the same guarantees as computations with floating-point numbers. So, in our opinion, Chebfun is a good role model for neural operators. In the present work, we use only a small part of Chebfun functionality, possible extensions are discussed in \cref{section:Conclusion_and_future_work}.

\section{Consequences of spectral representation}
\label{section:Consequences_of_spectral_representation}
Sufficiently smooth functions are indistinguishable from their Fourier and Chebyshev series if a reasonable number of terms is used. Given that, we restrict our attention to finite series and use duality between samples and coefficients to illustrate the accumulation of aliasing errors in FNO and discuss super-resolution. In what follows we mainly consider Fourier series and refer to $\phi(x; k) = \exp(i \pi k x)$ as harmonic with frequency $k$.

\subsection{Activation functions and aliasing}
\label{subsection:Activation_functions_and_aliasing}
Suppose our function $f(x)$ is (exactly) represented as Fourier series with $|k| < N$ terms. We can equivalently store values of the function on the uniform grid with $2N + 1$ points. However, when we apply activation function $\sigma(x)$, the resulting function $\sigma(f(x))$ may contain higher frequencies. If this is the case, composition $\sigma(f(x))$ can not be represented on the same grid with $2N+1$ points. The higher harmonics with $k>N$ can be computed on grid with $2N + 1$ points, but it is known that they are indistinguishable from lower harmonics with $k<N$ (see \cite[Chapter 11, Definition 23]{boyd2001chebyshev} for more details about Fourier aliasing and \cite[Chapter 4, Theorem 4.1]{trefethen2019approximation} for Chebyshev aliasing). If all operations are poinwise as in case of DeepONet or PiNN, it is not problematic. However, for FNO which uses global information, aliasing will present a problem. To quantify distortions introduced by activation function we define relative aliasing error $E_{a}$ for the Fourier series\footnote{For Chebyshev series norm is induced by scalar product with weight $1\big/\sqrt{1-x^2}$.}
\begin{equation}
    \label{aliasing_error}
    \sigma(f(x)) \equiv \sum_{i} c_i \phi_i,~E_a(f, \sigma)\equiv\frac{\left\|\sum_{|i|>N}c_i \phi_i\right\|_2}{\left\|\sigma(f(x))\right\|_2}.
\end{equation}
Aliasing error defined in \cref{aliasing_error} measures the norm of harmonics we cannot possibly resolve on the given grid relative to the norm of the function, transformed by activation. The following result gives aliasing error for rectifier and two extreme basis functions.

\begin{figure*}
    \centering
    \subfloat[]{
        \label{fig:FNO_aliasing}
        \includegraphics[width=0.48\textwidth]{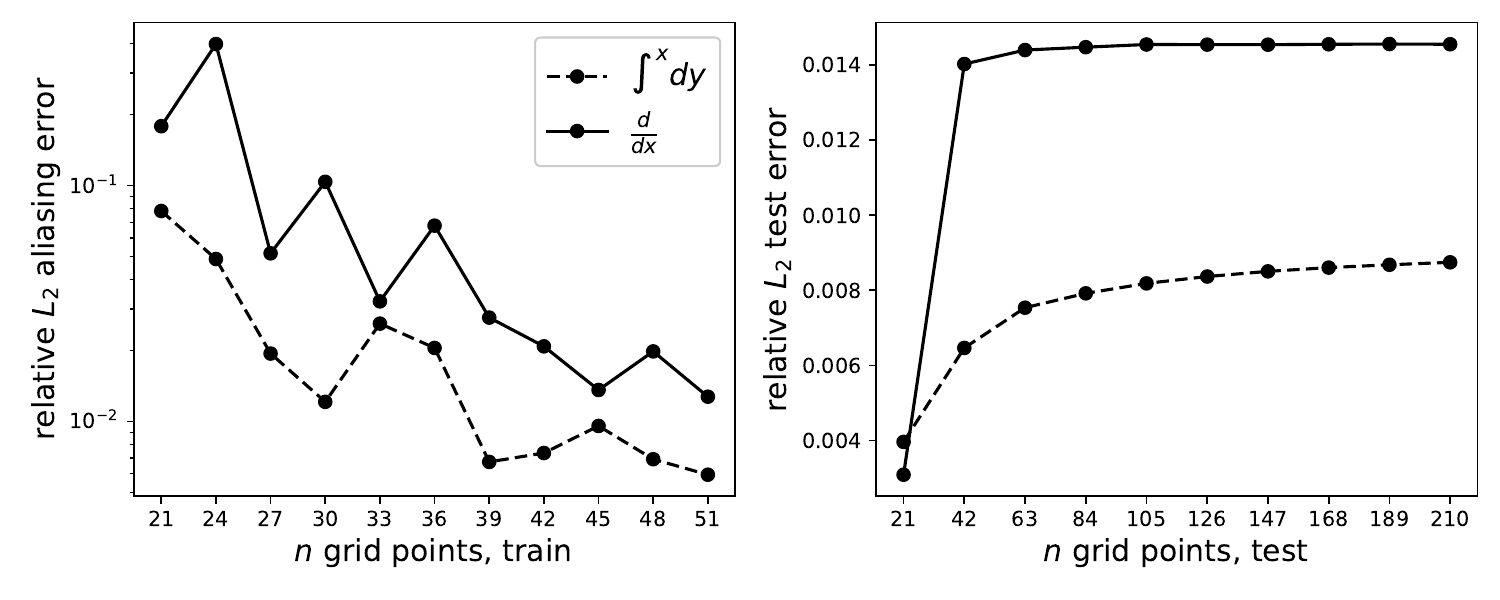}
    }
    \subfloat[]{
        \label{fig:super_resolution}
        \raisebox{-2pt}{\includegraphics[width=0.48\textwidth]{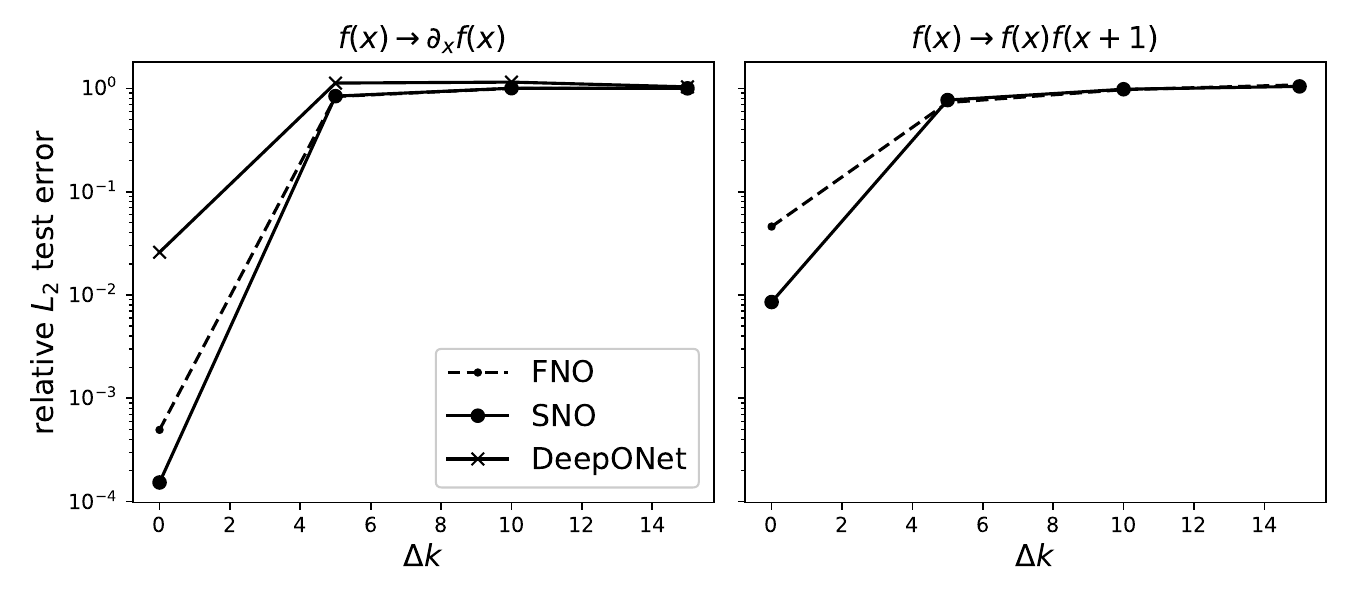}}
    }
    \caption{(a) On the left, relative aliasing error for FNO $N(u, x)$, computed as $\left\|N(u(x_{2h}), x_{2h}) - \left[N(u(x_{h}), x_{h})\right]_{2h}\right\|$ normalized on $\left\|\left[N(u(x_{h}), x_{h})\right]_{2h}\right\|$, where $\left[\cdot\right]_{2h}$ is a projection on grid with spacing $2h$. FNO was trained on the grid with $n$ point ($x$ axis) on functions with frequencies $[0, 10]$. Aliasing error decreases when larger grids are used. On the right, relative test error for the input with $n$ point ($x$ axis) for FNO trained for input with $n=21$ points. Relative test errors is larger for finer grid, which means FNO has systematic bias caused by aliasing. (b) Neural operators were trained on input with harmonics $[k_{\min}, k_{\max}]$ and grid with $100$ points and evaluated on input with harmonics $[k_{\min} + \Delta k, k_{\max} + \Delta k]$ on grid with $300$ points. Relative test error sharply increases with $\Delta k$, super-resolution is not observed.}
    \label{fig:super_resolution_aliasing}
\end{figure*}

\begin{theorem}[Aliasing and ${\sf ReLU}$]
\label{theorem:aliasing_and_relu}
On the uniform grid with $2N + 1$ points and on Chebyshev grid with $N+1$ points for $E_a$ defined in \cref{aliasing_error} we have
\begin{equation}
    E_a(\cos(\pi N x), {\sf ReLU}) = E_a(T_N(x), {\sf ReLU}) = \frac{1}{\pi}\sqrt{\frac{\pi^2}{2} - 4}\simeq 0.31
\end{equation}
\end{theorem}
The proof can be found in \cref{appendix:aliasing_and_relu}. More results on aliasing for composition with smooth functions can be found in \cite{bergner2006spectral}. The aliasing error is quite substantial, but since all energy in the theorem above is confined in highest possible harmonics, in practice one can expect to have milder discrepancy.

\cref{theorem:aliasing_and_relu} predicts that in certain situations FNO has a systematic bias. Namely, if one trains FNO on an insufficiently fine grid, activation functions introduce distortions that FNO will learn to mitigate. When the grid is sufficiently refined, aliasing errors disappear, but since FNO was trained to mitigate them, it predicts output functions that differ from targets it was trained on. The precise discrepancy is hard to analyze but one can expect that for operators with smoother output aliasing will be milder. This hypothesis is confirmed by our experiments. Indeed, the right graph in \cref{fig:FNO_aliasing} shows that for the integration operator, which performs smoothing, aliasing leads to an approximately two-fold error increase, while for the differentiation operator the increase is five-fold. For both examples, FNO still produces good approximations, but one should be cautious using FNO for more complex problems, because in these applications FNO tends to be less accurate \cite{pathak2022fourcastnet}, \cite{grady2022towards}.

Our solution to the problem of aliasing errors is to decouple interpolation from the mapping between functional spaces. It is described in \cref{section:Spectral_neural_operator}. However, it is possible to decrease aliasing error even when FNO is used for interpolation. One simply needs to train (or fine-tune) on a sufficiently fine grid. The decrease of aliasing error in this scenario is illustrated in the left plot of \cref{fig:FNO_aliasing}.

\subsection{Super-resolution}
\label{subsection:Superresolution}

In a series of papers \cite{li2020fourier}, \cite{kovachki2021neural}, \cite{li2020neural} authors claim that neural operators can handle different resolution, and even demonstrate ``zero-shot super-resolution''. The definition of the later was not given, but ``zero-shot super-resolution'' is tested in two steps. First, FNO is trained on data on the uniform grid with $N$ points and reaches test error $E^{N}_{\text{test}}$. After that FNO is evaluated on the same data but on grid with $M>N$ points and demonstrates test error $E^{M}_{\text{test}} \simeq E^{N}_{\text{test}}$. This procedure is problematic for two reasons.

First, high-resolution data is used as an input during the second stage. It is abnormal, because in classical works on super-resolution several low-resolution images are combined to improve resolution \cite{tian2011survey} without the direct use of high-resolution data. Second, the presented procedure implies that resolution increases when the number of grid points increases. This is obviously not the case in general. To claim such a thing, one needs to specify the family of function she is working with. For example, consider functions given by trigonometric series with $2N+1$ terms. If we sample these functions on the uniform grid with $2N+1$ points, we can exactly recover all $2N+1$ coefficients and no further refinement of the grid is going to provide novel information about the function, i.e., to increase resolution.

Given that, we provide an improved super-resolution test that also measures the generalization properties of neural operators. At the first stage one trains neural operator on functions with harmonics $|k|\leq K$ on the grid with $N \geq 2K + 1$ points.\footnote{Here we refer to complex exponential $\phi_{k}(x)$ and uniform grid. The extension to the case of Chebyshev polynomials is obvious.} At the second stage, one evaluates neural operator on functions with some or all harmonics $|k| > K$ on the grid with $M \geq (2\max|k| + 1)$ points. Note, that our test does not take into account the properties of the operator to be learned. In some cases the test can produce meaningless results, see \cref{appendix:Additional_comments_on_super_resolution} for details.

\begin{wraptable}{r}{0.25\textwidth}
    \caption{Neural operators trained on high harmonics and evaluated on low harmonics.}
    \input{tables/low_frequencies}
    \label{table:low_frequencies}
\end{wraptable}
We perform this test on FNO, DeepONet and spectral neural operator (see \cref{section:Spectral_neural_operator}) trained to take derivative and perform transformation $f(x)\rightarrow f(x)f(x+1)$. The results and details of the experiments are in \cref{fig:super_resolution}.\footnote{DeepONet fails to learn the second operator, that is why data for DeepONet is absent in the second picture.} Clearly, no super-resolution is observed for neural operators. This is a frustrating experience because very simple integration and differentiation rules can pass this test effortlessly. The other consequence of experiment \cref{fig:super_resolution} is that all neural operators learn only frequencies that are present during training. This is true for low frequencies too. To illustrate this, we train FNO and spectral neural operator to take indefinite integral, using harmonics $k\in[15, 25]$ and evaluate both neural operators on data with lower harmonics $k\in[15-\Delta k, 25-\Delta k]$. The results are given in \cref{table:low_frequencies}. Again, we can see that neural operators are unable to perform well on low harmonics when they are not presented during the training stage.

\section{Spectral neural operator}
\label{section:Spectral_neural_operator}
Both results in \cref{table:low_frequencies} and in \cref{fig:super_resolution} suggest that neural operators perform mapping of finite-banded functions.\footnote{One may object that we considered only smooth data. This is, indeed, true, but as shown in \cite{de2022cost} both FNO and DeepONet produce Runge-like oscillations when applied to data with discontinuities, so with current versions of these operators non-smooth data is out of question.} Bands are fixed by the choice of architecture and training data. Then, the natural direction is to define neural operator with fixed number of harmonics (basis functions) for both input  and output function:
\begin{equation}
    \label{sno}
    \sum_{i}g_{i}(x) d_i = f_{\sf in}(x) \overset{\sf N}{\longrightarrow} f_{\sf out}(x) = \sum_{i}g_{i}(x) b_i,
\end{equation}
where $g_{i}(x)$ are either Chebyshev polynomial or complex exponential. So, neural network maps finite number of input coefficients $\left\{d_i\right\}$ to the finite number of output coefficients $\left\{b_i\right\}$. Using the duality between samples on the grid (Chebyshev or uniform) and values of coefficients we can alternatively map on-grid samples. We coin term {\em Spectral Neural Operator} (SNO) for any network that provide mapping specified in \cref{sno} and preserve the structure of the series.

It is important to understand \cref{sno} implies that SNO only realizes mapping between two functions given by truncated series. If one needs to compute these functions on the finer grid, Chebyshev or trigonometric interpolation should be use. In the same vein, Chebyshev/trigonometric smoothing (the chopping of the series) is used on the input function, if it is not conformable with a particular architecture of SNO.

\input{Tikz_pictures/Chebyshev_LA}

The other important point is that series from \cref{sno} are truncated, so coefficients $d_i$ and $b_i$ can be stacked in finite-dimensional vectors. However, since these vectors represent sequences, linear algebra becomes slightly unusual. For example, it is possible to add two vectors even if their length does not coincide, or to multiply a vector by a matrix not conformable with it. The simplest way to understand this behavior is to think about sequences as vectors with an infinite number of trailing zeros that represent (zero) coefficients of truncated basis functions. They also can be considered as a simplified version of quasimatrices \cite{townsend2015continuous}. Standard linear algebra operations for sequences are illustrated in \cref{table:basic_operations_1}. It is also explained there, that, in the space of coefficients, linear operators quite appropriately becomes low-rank integral operators.

Except for unusual linear algebra, coefficients in \cref{sno} can be treated as ordinary vectors, so essentially any network that can map $\mathbb{C}^{n}$ to $\mathbb{C}^{m}$ qualifies as a neural operator.\footnote{The mapping can be between coefficients or, thanks to duality, between values of the function on the Chebyshev/uniform grid. When the later mapping is chosen, interpretation given in \cref{table:basic_operations_1} does not apply.} In the next section we gather benchmarks for SNO realized as ordinary feedforward neural networks, FNO and DeepONet.

\section{Numerical examples}
\label{section:Numerical_examples}

\begin{table*}[t]
    \centering
    \caption{Relative $L_2$ test errors for various problem. SNO(Ch) -- Spectral Neural Operator in Chebyshev basis, SNO(F) -- Spectral Neural Operator in the Fourier basis (in case of $2D$ problems, denoted by $xt$, Chebyshev basis is used along $t$ if the problem is not periodic), FNO \cite{li2020fourier} -- Fourier Neural Operator, DeepONet \cite{lu2019deeponet} -- Deep Operator Network.}
    \input{tables/truncated_big_table_with_all_results}
    \label{tab:big_table_with_all_results}
\end{table*}

Here we compare FNO \cite{li2020fourier}, DeepONet \cite{lu2019deeponet}, SNO(Ch) given by \cref{sno} in Chebyshev basis, SNO(F) given by \cref{sno} in Fourier basis. Precise architectures and details on training are described in \cref{appendix:details_on_training_and_architectures}. In addition, we perform benchmarks with more elaborate SNO architectures in \cref{appendix:benchmarks_for_additional_architectures}, and analyze sensitivity to initialization in \cref{appendix:sensitivity_to_initialization}. Code is available in \href{https://github.com/VLSF/SNO}{https://github.com/VLSF/SNO}.

In most of our benchmarks we used the following family of random functions
\begin{equation}
    \label{random_functions}
    \begin{split}
    &R(k_{\min}, k_{\max}, \sigma) = \left\{\frac{\mathcal{R}\left(\sum_{k=k_{\min}}^{k_{\max}}d_{k}^{j}e^{i\pi k x}\right)}{\sqrt{\sum_{k=k_{\min}}^{k_{\max}}\left|d_{k}^{j}\right|^{2}}}\right\}_{j},\,\mathcal{R}(d^{j}),\,\mathcal{I}(d^{j})\sim \mathcal{N}(0, \sigma^2),\,j\in[1, N]
    \end{split}
\end{equation}
where $\mathcal{R}(\cdot)$, $\mathcal{I}(\cdot)$ are real and imaginary part of complex number, $\mathcal{N}(a, b)$ is a normal distribution with mean $a$ and variance $b$. Operators that we used for training are described below.

\subsection{Operators}
\subsubsection{Integration, shift and differentiation}
\label{section:operators.subsection:integration_shift_and_differentiation}
First three test operators are
\begin{equation}
    \begin{split}
        &f(x) \rightarrow \int^{x} f(y) dy,\,f(y)\in R(1, 10, 2);\,f(x) \rightarrow f(x)f(x+1),\, f(x)\in R(0, 15, 2);\\
        &f(x) \rightarrow \frac{d}{dx}f(x),\, f(x)\in R(0, 10, k_{\max}),\,k_{\max}\in\left\{10, 20\right\}.
    \end{split}
\end{equation}

\subsubsection{Parametric ODE}
\label{section:operators.subsection:parametric_ode}
In this case we are interested in finding mapping $f(t)\rightarrow y(t)$, i.e., the solution to parametric family of ODEs (modification of DETEST problem $A4$ \cite{hull1972comparing})
\begin{equation}
    \frac{d y(t)}{dt} = y(t)f(t),\,y(0) = \exp\left(\int^{t=0} f(\tau) d\tau\right),\,f(t)\in R(1, 30, 2).
\end{equation}

\subsubsection{Elliptic equation}
\label{section:operators.subsection:elliptic_equation}
Here we considered two elliptic equations. In both cases, the goal is to predict the solution from the diffusion coefficient $k(x)$. In $1D$ we have
\begin{equation}
    \label{elliptic_1D}
    -\frac{d}{dx} k(x) \frac{d}{dx} u(x) = 1,\,u(-1)=u(1)=0,\,k(x) = 10(\tanh(g(x))+1) + 1,\,g(x)\in R(0, 20, 2).
\end{equation}
Diffusion coefficient changes from $1$ to $20$ but smoothly, so this equation is a variation of the benchmark (Darcy flow) used in \cite{li2020fourier}. In $2D$ we consider
\begin{equation}
    \label{elliptic_2D}
    \begin{split}
    &-\left(\frac{d}{dx} k(x)k(y) \frac{d}{dx} +\frac{d}{dy} k(x)k(y) \frac{d}{dy} \right)u(x, y) = 1,\,k(x) = \left(3(\tanh(g(x))+1) + 1\right),
    \end{split}
\end{equation}
with zero Dirichlet boundary conditions inside the square $[-1, 1]^2$ and $g(x),\,g(y)$ from $R(0, 20, 2)$.

\subsubsection{Burgers equation}
\label{section:operators.subsection:burgers_equation}
The goal is to infer solution of $ \frac{\partial}{\partial t} u(x, t) + \frac{1}{2}\frac{\partial}{\partial x}u(x, t)^2 - \nu \frac{\partial^2}{\partial x^2}u(x, t) = 0,\,u(x, 0) = f(x),\,x\in[-1,1]$, $f(x) \in R(0, 20, 2)$. In $1D$ we learn mapping $f(x)\rightarrow u(x, t=1)$, in $2D$ -- mapping $f(x)\rightarrow u(x, t),\,t\in[0, 1],\,x\in[-1, 1]$. We consider less viscous $\nu=0.01$ and more viscous $\nu=0.1$ scenarios.

\subsubsection{KdV equation}
\label{section:operators.subsection:kdv_equation}
In this benchmark we consider KdV equation
$\frac{\partial}{\partial t} u(x, t) + p\frac{\partial}{\partial x}u(x, t)^2 + \frac{\partial^3}{\partial x^3}u(x, t),\,u(x, 0) = f(x).$
More precisely, we use two exact solutions:
\begin{enumerate}
    \item Propagation of the single soliton ($p=1$) \cite{fornberg1978numerical}:
        \begin{equation}
            \label{soliton}
            u(x, t) = 3\left(\frac{a}{\cosh\left(a(x+x_0)\big/2 - a^3 t\big/2\right)}\right)^2,
        \end{equation}
    here $a$ defines amplitude and traveling speed and $x_0$ defines location of the maximum at $t=0$. We sample $a$ from uniform distribution on $[10, 25]$ and $x_0$ -- from uniform distribution on $[-1, 1]$. The goal is to predict $u(x, t=0.001)$ from $u(x, t=0)$ in $1D$ case, and $u(x, t),\,t\in[0,0.001]$ from $u(x, t=0)$ in $2D$ case.
    \item Interaction of two solitons ($p=6$) \cite{yoneyama1984korteweg}:
        \begin{equation}
            \label{two_solitons}
            \begin{split}
                &\phi_{i} = a_{i}(x+x_0^{i}) - 4a^3_{i} t,\\
                &\psi_1(x, t) = 2a_2^2(a_1^2 - a_2^2)\sinh\phi_1,\,
                \psi_2(x, t) = 2a_1^2(a_1^2 - a_2^2)\cosh\phi_2,\\
                &u(x, t) = \frac{\psi_1(x, t) + \psi_{2}(x, t)}{\left(a_1\cosh\phi_1\cosh\phi_2 - a_2\sinh\phi_1\sinh\phi_2\right)^2},\\
            \end{split}
        \end{equation}
        here $a_{i}$ defines the speed and amplitudes of solitons and $x^{i}_0$ their initial locations. We sample $a_1$ from uniform distribution on $[10, 25]$, $a_2=\chi a_1$, where $\chi$ is from uniform distribution on $[1/2, 1]$, the initial location was fixed to $x_0^{1}=-0.6$, $x_{0}^{2}=-0.5$. This way faster soliton is always behind the slower one and the neural operator needs to predict the whole interaction dynamics. For this problem time interval is $[0, 0.005]$, so the task is, again, to predict the solution in the final step, or on the whole spatio-temporal domain.
\end{enumerate}

\subsubsection{Nonlinear Schr\"odinger}
\label{section:operators.subsection:nonlinear_schrodinger}
Here we consider Kuznetsov-Ma breather \cite{kedziora2012second} given by
\begin{equation}
    \label{breather}
    \psi(x, t) = \left(\frac{-p^2\cos(\omega t) -2 i p\nu \sin(\omega t)}{2\cos(\omega t) - 2\nu \cosh(px\big/\sqrt{2})} - 1\right)e^{it},\,p=2\sqrt{\nu^2-1},\,\omega=p\nu.
\end{equation}
The goal is to predict $\left\|\psi(x, t)\right\|$ on $x\in[-1, 1]$, $t\in[0, 5]$ from initial conditions. The single parameter $\nu$ is sampled from uniform distribution on $[1.5, 3.5]$.

\subsection{Discussion of results}
\label{section:operators.subsection:discussion_of_results}
All results are gathered in \cref{tab:big_table_with_all_results}.

As we can see for many problems SNO performs better or comparable to FNO. The only exception is viscose Burgers equation with low viscosity. On this particular example both SNO networks failed to achieve reasonable test error. DeepONet perform much worse on our benchmarks. Especially, when high frequencies are present in the target. However note, that we use basic version of DeepONet. When more elaborate versions of DeepONet are used train error becomes comparable with the one of FNO \cite{lu2022comprehensive}.

\section{Conclusion and future work}
\label{section:Conclusion_and_future_work}
Using basic spectral techniques we describe the aliasing bias of FNO, discuss generalization properties and super-resolution of neural operators, and define SNO -- neural operator that operates on sequences of Chebyshev or Fourier coefficients. Benchmarks show that SNO is competitive when compared with other neural operators. So SNO appears to be a promising candidate for further investigations.

Limitations and natural lines for improvement are listed below.
\begin{enumerate}
    \item SNO in its current version is only applicable to smooth input and output data, because of the Gibbs phenomenon.\footnote{The same is true for DeepONet and FNO \cite{de2022cost}.} It is known that the Gibbs phenomenon can be completely resolved for spectral methods with a special reconstruction procedure \cite{gottlieb1992gibbs}. We plan to study this and related techniques in application to neural operators.
    \item Basis functions that we used are completely non-adaptive. In further research we plan to study how to implement an adaptive global basis using butterfly matrices \cite{dao2019learning}, which are known to be connected to fast special function transforms \cite{o2010algorithm}.
    \item The architecture that we used is very basic. It does not exploit in full the rich set of operations described in \cref{section:Chebyshev_and_trigonometric_polynomials}. We plan to implement and test substantially more complex SNO, possibly with adaptive operations defined in Chebfun\footnote{\href{https://www.chebfun.org}{https://www.chebfun.org}} \cite{trefethen2007computing}.
\end{enumerate}

\printbibliography[
heading=bibintoc,
title={Bibliography}
]

\appendix
\section{Aliasing and ReLU}
\label{appendix:aliasing_and_relu}
Here we prove the result stated in \cref{theorem:aliasing_and_relu}.
\begin{proof}
It is easy to see that
\begin{equation}
    {\sf ReLU}(x) = \sum_{i=0}^{\infty} p_i T_i(x),\,p_i =
    \begin{cases}
    \frac{1}{\pi},\,i=0;\\
    \frac{1}{2},\,i=1;\\
    \frac{2}{\pi}\frac{\cos(\pi i /2)}{1 - i^2},\,i\geq2,
    \end{cases}
\end{equation}
holds on the interval $[-1, 1]$. Since both $\cos(\pi N x)$ and $T_{N}$ have $[-1, 1]$ as their codomain, we can use Chebyshev series of ${\sf ReLU}$ to compute composition. For Chebyshev polynomial we use composition formula $T_{n}(T_{m}(x)) = T_{nm}(x)$ and for Fourier we use trigonometric definition $T_n(\cos(\theta)) = \cos(n\theta)$ to obtain
\begin{equation}
    \begin{split}
    &{\sf ReLU}(T_{N}(x)) = \sum_{i=0}^{\infty}p_i T_{iN}(x),\,{\sf ReLU}(\cos(\pi N x)) = \sum_{i=0}^{\infty}p_i \cos(\pi iN x).
    \end{split}
\end{equation}
Using norm of Chebyshev polynomial (with weight $1\big/\sqrt{1-x^2}$) $\left\|T_{i}(x)\right\|^2 = \pi \big/ (2 - \delta_{i0})$ and norm of trigonometric polynomial $\left\|\cos(\pi n x)\right\|^2_2 = (1 + \delta_{i0})$ we obtain for both aliasing errors
\begin{equation}
    \label{aliasing_error_p}
    E_{a}(\cdot, {\sf ReLU}) = \sqrt{\frac{\sum_{i=2}^{\infty} p_i^2}{2p_0^2 + p_1^2 + \sum_{i=2}^{\infty} p_i^2}}.
\end{equation}
To compute infinite sum $\sum_{i=2}^{\infty} p_i^2$ we can use Parseval's identity
\begin{equation}
    2p_0^2 + p_1^2 + \sum_{i=2}^{\infty} p_i^2 = \int \left({\sf ReLU}(\cos(\pi N x))\right)^2 dx = \frac{1}{2} \int \left(\cos(\pi N x)\right)^2 dx,
\end{equation}
where the last equality follows from simple geometric considerations. From Parseval's identity we find $\sum_{i=2}^{\infty} p_i^2=(\pi^2 - 8)\big/4\pi^2$ and confirm the result.
\end{proof}

It is also possible to extend results of \cref{theorem:aliasing_and_relu} considering refinement of the grid, i.e., aliasing errors of the form
\begin{equation}
    \sigma(f(x)) \equiv \sum_{i} c_i \phi_i,~E_a(f, \sigma, k)\equiv\frac{\left\|\sum_{|i|>kN}c_i \phi_i\right\|_2}{\left\|\sigma(f(x))\right\|_2},\,k\geq1.
\end{equation}
From practical standpoint $E_a(f, \sigma, k)$ measure aliasing error if one first performs (spectral) interpolation of $f(x)$ on grid with $kN$ points and after that applies activation function $\sigma$. To compute $E_a(f, \sigma, k)$ one needs to replace the sum $\sum_{i=2}^{\infty} p_i^2$ in the numerator of \cref{aliasing_error_p} with $\sum_{i=(k+1)}^{\infty} p_i^2$.

\section{Additional comments on super-resolution}
\label{appendix:Additional_comments_on_super_resolution}
An informal definition of super-resolution from \cref{subsection:Superresolution} combined with relative $L_2$ error used to measure performance can lead to erroneous conclusions in situations when high-harmonics do not matter.

\begin{wraptable}{r}{0.5\textwidth}
    \caption{Neural operators trained on harmonics $k\in[1, 30]$ and tested on $k\in[1 + \Delta k, 30 + \Delta k]$.}
    \input{tables/super_resolution_additional}
    \label{table:superresolution_additional}
\end{wraptable}
To exemplify, we perform the same super-resolution test as in \cref{subsection:Superresolution} form problem \cref{section:operators.subsection:parametric_ode}. According to the results, given in \cref{table:superresolution_additional}, neural operators appear to show super-resolution. All of them slightly deteriorate with the increase of $\Delta k$ but overall the performance remains good. To understand the reason behind that, we analyze the behavior of the solution.

For $f(t) = \mathcal{R}\left(\sum_{j=1 + \Delta k}^{30 + \Delta k} c_{j} \exp(ij\pi t)\right)$ exact solution to ODE from \cref{section:operators.subsection:parametric_ode} reads
\begin{equation}
    y(t) = \exp\left(\mathcal{R}\left(\sum_{j=1 + \Delta k}^{30 + \Delta k} \frac{c_{j}}{ij\pi} \exp(ij\pi t)\right)\right).
\end{equation}
As $\Delta k$ increases, the absolute value of the argument of exponent decreases, so the solution is well approximated by
\begin{equation}
    y(t) \simeq 1 + \mathcal{R}\left(\sum_{j=1 + \Delta k}^{30 + \Delta k} \frac{c_{j}}{ij\pi} \exp(ij\pi t)\right).
\end{equation}
Since $|c_j| \simeq 1$, it is enough for neural operator to predict constant output $1$ to have low relative $L_2$ error. In this sense the family of ODEs from \cref{section:operators.subsection:parametric_ode} simply have trivial high-frequency dynamics. The same is true for elliptic problem from \cref{section:operators.subsection:elliptic_equation}. It should be possible to extend the super-resolution test for these problems. But the analysis becomes cumbersome and, probably, irrelevant for applications.

\section{Details on training and architectures}
\label{appendix:details_on_training_and_architectures}
\subsection{Architectures}
Here we detail architectures used in the main text of the article as well as in \cref{appendix:benchmarks_for_additional_architectures}. Code is available in \href{https://github.com/VLSF/SNO}{https://github.com/VLSF/SNO}.
\subsubsection{FNO}
Architecture for FNO is taken from \cite{li2020fourier}.

For $D=1$ problems FNO consists of linear layer that transforms $2$ input features ($x$ and input function $f(x)$) into $64$ features, then $4$ layers with integral operator follow. In each of these layers $16$ Fourier modes are used. Final linear layer transforms $64$ features into $1$ (output functions).

For $D=2$, the construction is essentially the same, but the number of features is $32$, and the number of modes is $12$.

FNO has ${\sf ReLU}$ as activation function.

\subsubsection{DeepONet}

Architecture for DeepONet is from \cite{lu2019deeponet}.

For $D=1$ problems stacked branch network has $4$ dense layers with $100$ neurons in each layer, trunk net also has $4$ layers with $100$ features.

In $2D$ stacked branch net contain $100$ neurons for each dimension and trunk net has $100$, $100$, $100$, $10^4$ features in four layers.

DeepONet has $\tanh$ as an activation function.

\subsubsection{SNO}

SNO operates on coefficients of Chebyshev or Fourier series. The input and output of typical layer is a few functions, so for convenience of explanation we form rectangular (tall) matrix stacking them together
\begin{equation}
    \label{matrix_with_functions}
    U =
    \begin{pmatrix}
        | & & | \\
        u_{1} & \dots & u_{n}\\
        | & & | \\
    \end{pmatrix} \in \mathbb{C}^{m \times n}.
\end{equation}

In \cref{matrix_with_functions} each column correspond to a function and row $j$ correspond to coefficient $c_{j}$ in either Chebyshev or Fourier series.

\subsubsubsection{SNO(Ch), SNO(F), xSNO(Ch), xSNO(F)}
The first type of architecture that we use resembles the design of FNO \cite{li2020fourier}, but with parameter sharing in feature space. The network consists of three feedforward neural networks
\begin{equation}
    U \overset{N_1}{\longrightarrow} V \overset{N_2}{\longrightarrow} Y \overset{N_3}{\longrightarrow} W,
\end{equation}
where $U$, $V$, $Y$, $W$ have the same structure as the matrix in \cref{matrix_with_functions}.

Neural networks $N_1$ and $N_3$ consists of layers
\begin{equation}
    \begin{split}
    &U^{(n+1)} = \sigma\left(U^{(n)}A + b\right)\in \mathbb{C}^{k\times m}\\
    &U^{(n)} \in \mathbb{C}^{k\times l},\,A \in \mathbb{C}^{l\times m},\,b\in\mathbb{C}^{1\times m},
    \end{split}
\end{equation}
where $\sigma$ is a pointwise activation. Note that $b$ is not conformable with $U^{(n)}A$ but rows represent sequences, so this operation is well-defined. It corresponds to the addition of constant functions. As we can see networks $N_1$, $N_3$ can only change the number of functions and apply pointwise operations to functions. As in case of FNO, network $N_1$ is used to increase the number of features, and $N_2$ is used to decrease it to the number required by the output.

Neural network $N_2$ is stacked from the following layers
\begin{equation}
    \begin{split}
    &U^{(n+1)} = \sigma\left(BU^{(n)}A + b\right)\in \mathbb{C}^{r\times m}\\
    &B \in \mathbb{C}^{r\times k},\,U^{(n)} \in \mathbb{C}^{k\times l},\,A \in \mathbb{C}^{l\times m},\,b\in\mathbb{C}^{r\times m},
    \end{split}
\end{equation}
that is, in additional to operations presented in $N_1$ and $N_3$, it also contains linear operator $B$ that corresponds to low-rank integral operator described in \cref{table:basic_operations_1}.

In case $D>1$ the structure of $N_1$ and $N_3$ does not change, and in $N_2$ integral operators are applied to each dimension separately.

We have several SNO:
\begin{enumerate}
    \item SNO(Ch) uses Chebyshev basis along each dimension;
    \item SNO(F) uses Chebyshev basis along temporal (if time is present), and Fourier basis along spatial (more general, in periodic) dimension.
    \item xSNO(Ch) uses values of the function evaluated on Chebyshev grid.
    \item xSNO(F) uses values on Chebyshev grid for non-periodic dimensions (temporal, as a rule), and values on uniform grid otherwise.
\end{enumerate}

For $D=1$ problems SNO(Ch) and xSNO(Ch) has three layers with integral operator $100$ neurons in each, $N_1$ and $N_3$ consist of a single layer with $20$ features. SNO(F) and xSNO(F) has the same architecture, but the number of neurons in integral operator is $51$.

In $D=2$ architectures of SNO(Ch), xSNO(Ch), SNO(F), xSNO(F) repeat the one for $D=1$, with the same number of features, $100$ neurons per non-periodic dimension, $51$ neurons for periodic dimension.

As activation functions we use softplus applied separately to real and complex (if present) part of the input. For complex-valued networks we also tried more specialized activation functions described in \cite[Section III]{scardapane2018complex}, namely complex cardioid \cite[Equation 18]{scardapane2018complex} and complex-valued ${\sf ReLU}$ \cite[Equation 16  ]{scardapane2018complex} and did not observe any noticeable improvement in comparison with softplus.

\begin{table*}
    \centering
    \caption{Relative $L_2$ test errors for $D=1$ operators; xSNO(Ch) -- Spectral Neural Operator on Chebyshev grid, SxNO(F) -- Spectral Neural Operator on the uniform grid, xcSNO(Ch) -- combination of xSNO(Ch) and SNO(Ch), xcSNO(F) combination of xSNO(F) and SNO(F).}
    \input{tables/additional_architectures_1D}
    \label{tab:big_table_with_additional_architectures_1D}
\end{table*}

\subsubsubsection{xcSNO(Ch), xcSNO(F)}

We also consider more elaborate architecture that combines operations on the grid with operations on coefficients.

The network structure is as follows
\begin{equation}
    U \overset{N_2}{\longrightarrow} V\overset{F}{\longrightarrow} V_{c} \overset{N_3}{\longrightarrow} Y_{c} \overset{F^{-1}}{\longrightarrow} Y \overset{N_2}{\longrightarrow} W,
\end{equation}
where $U$ is an input function evaluated on Chebyshev or uniform grid, $N_2$ is as described above, $F$ is the transformation from values to coefficients, $N_3(V_{c}) = N_2(V_{c}) + V_{c}$ and each network $N_2$ has different weights.

We use these architectures only for $D=1$ problems. Both $N_2$ networks have two layers with $100$ (for Chebyshev basis) or $51$ (for Fourier basis) neurons, and $N_3$ has three layers with the same number of neurons. The number of features is $20$.

\subsection{Initialization}
We draw each bias from standard normal distribution. Values of each tensor that is used as linear operator $\sum_{l}A_{ijkl} b_{l}$ are also drawn from standard normal distribution but in addition are multiplied by normalization factor $\left(\sum_{l} 1\right)^{-1}$.

\subsection{Details on datasets and training}
\label{appendix:details_on_training_and_architectures.subsection:details_on_datasets_and_training}
In all benchmarks we use $100$ points on uniform grid, $100$ Chebyshev coefficients and $100$ Fourier harmonics (i.e., $51$ complex numbers in $D=1$ and $100\times51$ complex numbers in $D=2$) along each dimension. In the case of SNO, we recompute train and test errors on the uniform grid for it to be consistent with the results of other networks, working with functions on uniform grids. Also, in all benchmarks both train and test sets contain $1000$ examples. Networks were trained using an Adam optimizer with a learning rate $0.001$ with exponential decay $0.5$ applied each $10000$ epochs in $1D$ and $5000$ in $2D$. Training is performed with large batch sizes: $1000$ in $1D$ and $100$ in $2D$. The number of epochs in $1D$ is $50000$ for each operator. In $2D$ computational load becomes extensive so we train each network for roughly the same amount of (wall-clock) time $\sim 6.7$ hours. In terms of the number of epochs this is $50000$ for SNO(Ch), xSNO(Ch), xSNO(F), $\simeq 21600$ for FNO, $\simeq 31500$ for SNO(F), $\simeq 81400$ for DeepONet.

All training was performed on a single GPU NVIDIA Tesla V100 SXM2 $16$ GB.

For problems in \cref{section:operators.subsection:integration_shift_and_differentiation}, \cref{section:operators.subsection:parametric_ode}, \cref{section:operators.subsection:kdv_equation}, \cref{section:operators.subsection:nonlinear_schrodinger}, \cref{section:operators.subsection:nonlinear_schrodinger} solutions are explicitly known, so generation of datasets is straightforward.

For \cref{section:operators.subsection:burgers_equation} the exact solution is also available, but rarely used in practice. We use fourth-order Runge-Kutta with time step $0.0001$, pseudospectral discretization with trick from \cite{fornberg1978numerical} (see also \cite[Problem 10.6]{trefethen2000spectral}). Number of Fourier modes used is $100$.

For \cref{section:operators.subsection:elliptic_equation} we use standard Chebyshev pseudospectral method, i.e., spectral differentiation matrix \cite[Chapter 6]{trefethen2000spectral}. Matrix inversion was done with $LU$ decomposition. Chebyshev grid with $100$ points was used along each dimension.

\begin{table*}[t]
    \centering
    \caption{Relative $L_2$ test errors for $D=2$ operators; xSNO(Ch) -- Spectral Neural Operator on Chebyshev grid, xSNO(F) -- Spectral Neural Operator on the uniform grid, FNO \cite{li2020fourier} -- Fourier Neural Operator, DeepONet \cite{lu2019deeponet} -- Deep Operator Network.}
    \input{tables/additional_architectures_2D}
    \label{tab:big_table_with_additional_architectures_2D}
\end{table*}

\section{Benchmarks for additional architectures}
\label{appendix:benchmarks_for_additional_architectures}

In main text we show results for SNO(Ch) and SNO(F), i.e., neural operators defined in the space of coefficients. In \cref{tab:big_table_with_additional_architectures_1D} and \cref{tab:big_table_with_additional_architectures_2D} one can find results for architectures that transform values of the function on Chebyshev or uniform grid.

For $D=1$ problems one can see the results for additional architectures are not radically different from SNO(Ch) and SNO(F). The results for the smooth Burgers equation, Elliptic equation, and parametric ODE are slightly better. For Burgers equation with $\nu=0.001$, performance is still poor.

For $D=2$ problems the situation is largely the same. This time xSNO(Ch) performs slightly better than FNO on the Elliptic equation.

\section{Sensitivity to initialization}
\label{appendix:sensitivity_to_initialization}
Here we provide data on sensitivity to initialization. To reduce computational footprint, all results are presented for three selected $D=1$ problems. For each problem, we trained $10$ randomly initialized networks. As we can see from \cref{tab:sensitivity to initialization} training is somewhat sensitive to weight initialization, but not sensitive enough to alter conclusion in \cref{section:operators.subsection:discussion_of_results}.

\begin{table*}[t]
    \centering
    \caption{Sensitivity to initialization for selected $D=1$ problems. Below one can find $\mu(E) \pm \sigma(E)$, where $\mu(E)$ -- mean relative test error, $\sigma(E)$ -- square root of variance. Averaging is performed over $10$ random initializations. Situation $\mu(E)-\sigma(E) < 0$ means that error distribution is skewed to the right.}
    \resizebox{\columnwidth}{!}{\input{tables/sensitivity_to_initialization}}
    \label{tab:sensitivity to initialization}
\end{table*}
\end{document}

%% file: Tikz_pictures/tikz.tex
\usepackage{tikz}

\definecolor{vector_color}{RGB}{204, 204, 255}
\definecolor{extend_color}{RGB}{229, 255, 204}
\definecolor{truncate_color}{RGB}{255, 204, 255}

%% file: tables/low_frequencies.tex
\begin{tabular}{ccc}
    \toprule
    &\multicolumn{2}{c}{$E_{\text{test}}$}\\
    \cmidrule(lr){2-3}
    $\Delta k$ & FNO & SNO \\
    \hline
    $0$ & $0.12$ & $0.021$ \\
    $5$ & $0.61$ & $0.66$ \\
    $10$ & $1.0$ & $1.1$ \\
    \bottomrule
\end{tabular}

%% file: Tikz_pictures/Chebyshev_LA.tex
\begin{figure}
    \centerline{\begin{tabular}{c|c}
        Operation on vectors & Operation on functions
        \\
        \hline
        \\
        \input{Tikz_pictures/representation} &
        \begin{tabular}{l} \parbox{6cm}{Truncated series of Chebyshev polynomials can be represented as stacked vectors:}\\
        $\psi(x) = \sum\limits_{i=0}^{4} \psi_{i}T_{i}(x),$ $\chi_j(x) = \sum\limits_{i=0}^{3} \chi_{ij} T_{i}(x),$ \\
        \parbox{6cm}{
        \begin{multline*}
            \left|\psi\right> = \begin{pmatrix} \psi_0 & \psi_1 & \psi_2 & \psi_3 & \psi_4\end{pmatrix}^{T}, \\\left<\psi\right| = \frac{\pi}{2}\begin{pmatrix} 2\psi_0 & \psi_1 & \psi_2 & \psi_3 & \psi_4\end{pmatrix}.
        \end{multline*}}
        \end{tabular}
        \\
        \hline
        \\
        \input{Tikz_pictures/addition} & \begin{tabular}{l} \parbox{6cm}{Addition of two functions is equivalent to addition of two vectors with zero padding of the short one:}
        \\
        $\chi(x) = \sum\limits_{i=0}^{3} \chi_{i} T_{i}(x)$, $\psi(x) = \sum\limits_{i=0}^{7} \psi_{i} T_{i}(x),$ \\
        \parbox{6cm}{
        \begin{multline*}
            \psi(x) + \chi(x) = \sum\limits_{i=0}^{3} \left(\psi_{i} + \chi_i\right) T_{i}(x) \\+\sum\limits_{i=4}^{7} \left(\psi_{i} + 0\right) T_{i}(x).
        \end{multline*}}
        \end{tabular}
        \\
        \hline
        \\
        \input{Tikz_pictures/projection} &
        \begin{tabular}{l}
        \parbox{6cm}{Chebyshev scalar product of two functions is equivalent to scalar product of vectors with truncation of the long one:}
        \\
        $\chi(x) = \sum\limits_{i=0}^{1}\chi_{i} T_{i}(x)$, $\psi(x) = \sum\limits_{i=0}^{5}\psi_i T_{i}(x),$\\
        \parbox{6cm}{
        \begin{multline*}
         \left<\chi\,\right|\left.\psi\right> \equiv\int\limits_{-1}^{+1} \frac{dx}{\sqrt{1-x^2}} \chi(x) \psi(x)  \\ = \pi\chi_{0}\psi_{0} + \frac{\pi}{2} \chi_{1}\psi_{1}.
        \end{multline*}}
        \end{tabular}
        \\
        \hline
        \\
         \input{Tikz_pictures/integral_operator} &
        \begin{tabular}{l}
        \parbox{6cm}{Low-rank integral operator is equivalent to multiplication by low-rank matrix:}
        \\
        $\chi_{j}(x) = \sum\limits_{i=0}^{3}\chi_{ij} T_{i}(x)$, $\phi_{j}(x) = \sum\limits_{i=0}^{5}\phi_{ij} T_{i}(x),$
        \\
        $\psi(x) = \sum\limits_{i=0}^{4}\psi_i T_{i}(x),$
        \\
        \parbox[b]{5.7cm}{
        \begin{multline*}
            \sum\limits_{i=1}^{2}\left| \phi_i\right>_{k}\left<\chi_{i}\,\right|\left.\psi\right> = \sum\limits_{i=1}^{2}\sum\limits_{j=0}^{3}\phi_{ki}\chi_{ij}\psi_{j}.
        \end{multline*}}
        \end{tabular}
        \\
    \end{tabular}}
    \caption{Linear algebra for Chebyshev polynomials.}
    \label{table:basic_operations_1}
\end{figure}

%% file: Tikz_pictures/representation.tex
\begin{tikzpicture}[baseline=-50]
    \filldraw[fill=vector_color, draw=black] (0, 0) rectangle (0.5, -2.5);
    \draw[step=0.5, black, thin] (0, 0) grid (0.5, -2.5);
    \node (1) at (0.25, -3.0) {$\left|\psi\right>$};
    \node (2) at (1.5, -0.25) {$\left|\psi\right>_{1} = \psi_{0}$};
    \node (2) at (1.5, -2.25) {$\left|\psi\right>_{5} = \psi_{4}$};
    \filldraw[fill=vector_color, draw=black] (0, 0) rectangle (0.5, -2.5);
    \draw[step=0.5, black, thin] (0, 0) grid (0.5, -2.5);
    \filldraw[fill=vector_color, draw=black] (3.0, 0) rectangle (4.5, -2);
    \draw[step=0.5, black, thin] (3.0, 0) grid (4.5, -2);
    \node (1) at (3.35, -2.5) {$\left|\chi\right>_1$};
    \node (1) at (4.35, -2.5) {$\left|\chi\right>_3$};
\end{tikzpicture}

%% file: Tikz_pictures/addition.tex
\begin{tikzpicture}[baseline=-70]
    \filldraw[fill=vector_color, draw=black] (0, 0) rectangle (0.5, -2.0);
    \draw[step=0.5, black, thin] (0, 0) grid (0.5, -2.0);
    \node (1) at (1.0, -1.0) {$+$};
    \filldraw[fill=vector_color, draw=black] (1.5, 0) rectangle (2.0, -4.0);
    \draw[step=0.5, black, thin] (1.5, 0) grid (2.0, -4.0);
    \node (2) at (2.5, -1.0) {$=$};
    \filldraw[fill=vector_color, draw=black] (3.0, 0) rectangle (3.5, -2.0);
    \filldraw[fill=extend_color, draw=black] (3.0, -2.0) rectangle (3.5, -4.0);
    \draw[step=0.5, black, thin] (3.0, 0) grid (3.5, -4.0);
    \node (3) at (4.0, -1.0) {$+$};
    \filldraw[fill=vector_color, draw=black] (4.5, 0) rectangle (5.0, -4.0);
    \draw[step=0.5, black, thin] (4.5, 0) grid (5.0, -4.0);
    \node (4) at (4.75, -4.5) {$\left|\psi\right>$};
    \node (5) at (3.25, -4.5) {$\left|\chi\right>$};
\end{tikzpicture}

%% file: Tikz_pictures/projection.tex
\begin{tikzpicture}[baseline=-65]
    \filldraw[fill=vector_color, draw=black] (0, 0) rectangle (1.0, -0.5);
    \draw[step=0.5, black, thin] (0, 0) grid (1.0, -0.5);
    \filldraw[fill=vector_color, draw=black] (1.5, 0) rectangle (2.0, -1.0);
    \filldraw[fill=truncate_color, draw=black] (1.5, -1.0) rectangle (2.0, -3.0);
    \draw[step=0.5, black, thin] (1.5, 0) grid (2.0, -3.0);
    \node (2) at (2.5, -0.25) {$=$};
    \filldraw[fill=vector_color, draw=black] (3.0, 0) rectangle (4.0, -0.5);
    \draw[step=0.5, black, thin] (3.0, 0) grid (4.0, -0.5);
    \filldraw[fill=vector_color, draw=black] (4.5, 0) rectangle (5.0, -1.0);
    \draw[step=0.5, black, thin] (4.5, 0) grid (5.0, -1.0);
    \node (4) at (4.75, -1.5) {$\left|\psi\right>$};
    \node (5) at (3.5, -1.0) {$\left<\chi\right|$};
\end{tikzpicture}

%% file: Tikz_pictures/integral_operator.tex
\begin{tikzpicture}[baseline=-60]
    \filldraw[fill=vector_color, draw=black] (0, 0) rectangle (1.0, -3.0);
    \draw[step=0.5, black, thin] (0, 0) grid (1.0, -3.0);
    \filldraw[fill=vector_color, draw=black] (1.5, 0) rectangle (3.0, -1.0);
    \draw[step=0.5, black, thin] (1.5, 0) grid (3.0, -1.0);
    \node (1) at (2.25, -1.5) {$\left<\chi_{i}\right|$};
    \filldraw[fill=vector_color, draw=black] (3.5, 0) rectangle (4.0, -1.5);
    \filldraw[fill=truncate_color, draw=black] (3.5, -1.5) rectangle (4.0, -2.5);
    \draw[step=0.5, black, thin] (3.5, 0) grid (4.0, -2.5);
    \node (2) at (3.75, -3.0) {$\left|\psi\right>$};
    \node (2) at (0.5, -3.5) {$\left|\phi_{i}\right>$};
\end{tikzpicture}

%% file: tables/truncated_big_table_with_all_results.tex
\begin{tabular}{cccccc}
        \toprule
        dataset & SNO(Ch) & SNO(F) & FNO & DeepONet\\
        \hline
        Integration & $3.5\times10^{-4}$ & $\boldsymbol{2.6\times10^{-4}}$ & $3.3\times10^{-2}$ & $4.7\times10^{-2}$ \\
        $\partial_x$, $k_{\max}=10$ & $\boldsymbol{1.5\times10^{-4}}$ & $1.3\times10^{-3}$ & $5.0\times10^{-4}$ & $2.5\times10^{-2}$ \\
        $\partial_x$, $k_{\max}=20$ & $3.5\times10^{-4}$ & $\boldsymbol{2.6\times10^{-4}}$ & $3.3\times10^{-2}$ & $4.7\times10^{-2}$ \\
        $f\rightarrow f(x)f(x+1)$ & $\boldsymbol{8.5\times10^{-3}}$ & $6.5\times10^{-2}$ & $4.6\times10^{-2}$ & $1.8$ \\
        parametric ODE & $\boldsymbol{9.7\times10^{-4}}$ & $2.2\times10^{-3}$ & $1.1\times10^{-3}$ & $2.3\times10^{-2}$ \\
        Elliptic \eqref{elliptic_1D} & $\boldsymbol{1.5\times10^{-2}}$ & --- & $2.3\times10^{-2}$ & $2.2\times10^{-1}$ \\
        Burgers$_{x}$, $\nu=0.01$ & $5.8\times10^{-1}$ & $3.1\times10^{-1}$ & $\boldsymbol{6.3\times10^{-2}}$ & $6.0\times10^{-1}$ \\
        Burgers$_{x}$, $\nu=0.1$ & $\boldsymbol{1.1\times10^{-2}}$ & $2.8\times10^{-2}$ & $1.5\times10^{-2}$ & $1.4\times10^{-1}$ \\
        KdV$_{x}$ \eqref{soliton} & $3.2\times10^{-2}$ & $2.3\times10^{-2}$ & $\boldsymbol{1.3\times10^{-3}}$ & $4.6\times10^{-1}$ \\
        KdV$_{x}$ \eqref{two_solitons} & $1.0\times10^{-1}$ & $\boldsymbol{8.0\times 10^{-2}}$ & $9.0\times 10^{-2}$ & $2.3\times 10^{-1}$\\
        Elliptic \eqref{elliptic_2D} & $3.0\times10^{-2}$ & --- & $\boldsymbol{2.7\times10^{-2}}$ & $9.3\times10^{-2}$ \\
        Burgers${}_{xt}$, $\nu=0.01$ & $1.9\times10^{-1}$ & $2.6\times10^{-1}$ & $\boldsymbol{4.1\times10^{-2}}$ & $5.8\times10^{-1}$\\
        Burgers${}_{xt}$, $\nu=0.1$ & $3.4\times10^{-2}$ & $9.4\times10^{-2}$ & $\boldsymbol{6.5\times10^{-3}}$ & $4.1\times10^{-1}$\\
        KdV${}_{xt}$ \eqref{soliton}& $\boldsymbol{5.7\times10^{-3}}$ & $5.8\times10^{-3}$ & $6.5\times10^{-3}$ & $7.1\times10^{-2}$\\
        KdV${}_{xt}$ \eqref{two_solitons}& $\boldsymbol{1.9\times10^{-2}}$ & $2.1\times10^{-2}$ & $6.2\times10^{-2}$ & $1.6\times10^{-1}$\\
        Breather \eqref{breather} & $\boldsymbol{1.4\times 10^{-2}}$ & $5.2\times 10^{-2}$ & $2.6\times 10^{-2}$ & $1.5\times 10^{-1}$ \\
        \bottomrule
    \end{tabular}

%% file: tables/super_resolution_additional.tex
\begin{tabular}{cccc}
    \toprule
    &\multicolumn{3}{c}{$E_{\text{test}}$}\\
    \cmidrule(lr){2-4}
    $\Delta k$ & FNO & SNO & DeepONet \\
    \hline
    $0$ & $1.2\times 10^{-3}$ & $9.7\times 10^{-4}$ & $2.3\times 10^{-2}$\\
    $5$ & $2.8\times 10^{-3}$ & $4.9\times 10^{-3}$ & $3.0\times 10^{-2}$\\
    $10$ & $5.1\times 10^{-3}$ & $7.4\times 10^{-3}$ & $3.4\times 10^{-2}$\\
    $15$ & $5.7\times 10^{-3}$ & $9.7\times 10^{-3}$ & $4.1\times 10^{-2}$\\
    $20$ & $6.1\times 10^{-3}$ & $1.1\times 10^{-2}$ & $5.0\times 10^{-2}$\\
    $25$ & $6.1\times 10^{-3}$ & $1.4\times 10^{-2}$ & $6.6\times 10^{-2}$\\
    \bottomrule
\end{tabular}

%% file: tables/additional_architectures_1D.tex
\begin{tabular}{cccccc}
        \toprule
        dataset & xSNO(Ch) & xcSNO(Ch) & xSNO(F) & xcSNO(F)\\
        \hline
        Integration & $\boldsymbol{1.2\times10^{-3}}$ & $2.2\times10^{-3}$ & $8.5\times10^{-3}$ & $2.5\times10^{-3}$ \\
        $\partial_x$, $k_{\max}=10$ & $6.7\times10^{-4}$ & $8.7\times10^{-4}$ & $\boldsymbol{6.3\times10^{-4}}$ & $9.6\times10^{-4}$ \\
        $\partial_x$, $k_{\max}=20$ & $\boldsymbol{5.8\times10^{-4}}$ & $8.9\times10^{-4}$ & $6.3\times10^{-4}$ & $7.6\times10^{-4}$ \\
        $f\rightarrow f(x)f(x+1)$ & $1.2\times10^{-2}$ & $9.0\times10^{-3}$ & $\boldsymbol{6.9\times10^{-3}}$ & $9.3\times10^{-3}$ \\
        parametric ODE & $\boldsymbol{1.5\times10^{-4}}$ & $6.0\times10^{-4}$ & $8.3\times10^{-3}$ & $5.4\times10^{-4}$ \\
        Elliptic \eqref{elliptic_1D} & $\boldsymbol{1.1\times10^{-2}}$ & $7.9\times10^{-2}$ & --- & --- \\
        Burgers$_{x}$, $\nu=0.01$ & $\boldsymbol{4.1\times10^{-1}}$ & $5.8\times10^{-1}$ & $7.4\times10^{-1}$ & $4.2\times10^{-1}$ \\
        Burgers$_{x}$, $\nu=0.1$ & $\boldsymbol{2.4\times10^{-2}}$ & $2.5\times10^{-2}$ & $2.5\times10^{-2}$ & $4.0\times10^{-2}$ \\
        KdV$_{x}$ \eqref{soliton} & $6.7\times10^{-2}$ & $5.0\times10^{-2}$ & $3.6\times10^{-2}$ & $\boldsymbol{2.6\times10^{-2}}$ \\
        KdV$_{x}$ \eqref{two_solitons} & $8.2\times10^{-2}$ & $1.5\times10^{-1}$ & $1.1\times10^{-1}$ & $\boldsymbol{7.6\times10^{-2}}$\\
        \bottomrule
    \end{tabular}

%% file: tables/additional_architectures_2D.tex
\begin{tabular}{cccccc}
        \toprule
        dataset & xSNO(Ch) & xSNO(F) & FNO & DeepONet\\
        \hline
        Elliptic \eqref{elliptic_2D} & $\boldsymbol{2.2\times10^{-2}}$ & --- & $2.7\times10^{-2}$ & $9.3\times10^{-2}$ \\
        Burgers${}_{xt}$, $\nu=0.01$ & $1.5\times10^{-1}$ & $1.9\times10^{-1}$ & $\boldsymbol{4.1\times10^{-2}}$ & $5.8\times10^{-1}$\\
        Burgers${}_{xt}$, $\nu=0.1$ & $1.4\times10^{-2}$ & $1.8\times10^{-2}$ & $\boldsymbol{6.5\times10^{-3}}$ & $4.1\times10^{-1}$\\
        KdV${}_{xt}$ \eqref{soliton}& $1.9\times 10^{-2}$ & $1.3\times10^{-2}$ & $\boldsymbol{6.5\times10^{-3}}$ & $7.1\times10^{-2}$\\
        KdV${}_{xt}$ \eqref{two_solitons}& $\boldsymbol{2.7\times10^{-2}}$ & $3.1\times10^{-2}$ & $6.2\times10^{-2}$ & $1.6\times10^{-1}$\\
        Breather \eqref{breather} & $\boldsymbol{8.4\times 10^{-3}}$ & $9.0\times 10^{-3}$ & $2.6\times 10^{-2}$ & $1.5\times 10^{-1}$ \\
        \bottomrule
    \end{tabular}

%% file: tables/sensitivity_to_initialization.tex
\begin{tabular}{cccccc}
        \toprule
        dataset & SNO(Ch) & SNO(F) & FNO & DeepONet\\
        \hline
        Integration & $(8.6\pm3.7)\times10^{-4}$ & $(6.9\pm4.6)\times10^{-4}$ & $(3.5\pm1.6)\times10^{-2}$ & $(6.4\pm4.0)\times10^{-2}$ \\
        $\partial_x$, $k_{\max}=10$ & $(8.6\pm4.9)\times10^{-4}$ & $(1.3\pm0.6)\times10^{-3}$ & $(5.2\pm4.1)\times10^{-4}$ & $(2.5\pm0.1)\times10^{-2}$ \\
        parametric ODE &
        $(3.6\pm6.9)\times10^{-3}$ & $(2.2\pm0.5)\times10^{-3}$ & $(1.3\pm0.8)\times10^{-3}$ & $(3.2\pm0.5)\times10^{-2}$ \\
        \bottomrule
    \end{tabular}